\documentclass[12pt]{article}

\usepackage{amssymb,amsmath,latexsym,amsthm,epsfig,enumerate,graphics,hyperref}
\usepackage{setspace}
\usepackage{color}
\hypersetup{colorlinks=true}

\newcommand\ie{i.\,e.}

\def\SISPr{\mu}
\def\TruePr{\pi}
\def\coS{\overline{S}}

\def\F{\mathcal{F}}

\def\eps{\epsilon}

\newcommand{\X}{U}
\newcommand{\Y}{V}

\def\Prob#1{{\mathbf{Pr}\left({#1}\right)}}

\def\ProbCond#1#2{{\mathbf{Pr}\left({#1} \, {\big \vert} \, {#2} \right)}}

\def\ExpCond#1#2{{\mathbf{E}\left({#1} \mid {#2} \right)}}

\newtheorem{theorem}{Theorem}

\newtheorem{lemma}[theorem]{Lemma}
\newtheorem{note}[theorem]{Remark}

\newtheorem*{duplicate-theorem}{}

\newcommand{\T}{\mathbf{T}}
\newcommand{\Orc}{\Omega_{\mathbf{r},\mathbf{c}}}

\begin{document}

\date{June 24, 2011}

\author{
Ivona Bez\'{a}kov\'{a}\thanks{Department of Computer Science,
Rochester Institute of Technology, 102 Lomb Memorial Drive,
Rochester, NY 14623. Email: {\tt ib@cs.rit.edu}.} \and Alistair
Sinclair\thanks{Computer Science Division, University of California,
Berkeley, CA 94720. Email:~ {\tt sinclair@cs.berkeley.edu}.
Supported by NSF grants CCF-1016896 and CCF-0635153.} \and Daniel
\v{S}tefankovi\v{c}\thanks{Department of Computer Science,
University of Rochester, Rochester, NY 14627.  Email:~ {\tt
stefanko@cs.rochester.edu}.
Supported by NSF grant CCF-0910415.}
\and Eric Vigoda\thanks{College of
Computing, Georgia Institute of Technology, Atlanta, GA 30332.
Email:~ {\tt vigoda@cc.gatech.edu}.
Supported by NSF grants CCF-0830298 and CCF-0910584.}}

\title{\vspace*{-.2in}
Negative Examples for Sequential Importance \\ Sampling of Binary Contingency Tables}

\maketitle

\vspace*{-.4in}

\begin{center}
\\  Graphs with prescribed degree sequence;
Zero-one table
\end{center}

\begin{abstract}
The sequential importance sampling (SIS) algorithm
has gained considerable popularity for its empirical success.
One of its noted applications is to
the binary contingency tables problem, an important problem in statistics,
where the goal is  to estimate the number of $0/1$ matrices with prescribed row and column sums.
We give a family of examples in which the SIS
procedure, if run for any subexponential number of trials,
will underestimate the number of tables
by an exponential factor.
This result holds for any of the usual design choices in the SIS algorithm,
namely the ordering of the columns and rows.
 These are apparently the first theoretical results on the
efficiency of the SIS algorithm for binary contingency tables.
Finally, we present
experimental evidence that the SIS algorithm is efficient for row and column sums
that are regular.  Our work is a first step in determining the class of inputs
for which SIS is effective.
\end{abstract}

\section{Introduction}
Sequential importance sampling is a widely-used approach for
estimating the cardinality of a large set of combinatorial objects.
It has been applied in
a variety of fields, such as
protein folding
\cite{ZhangLiu2002},
population genetics
\cite{IGLRousset2005},
and signal processing
\cite{SigProc}.
Binary contingency tables is an application where the virtues of
sequential importance sampling have been especially
highlighted;
see Chen et al.~\cite{CDHL}.
This is the subject of this note.
Given a set of non-negative row sums
$\mathbf{r}=(r_1,\dots,r_m)$ and column sums $\mathbf{c}=(c_1,\dots,c_n)$,
let $\Omega=\Omega_{\mathbf{r},\mathbf{c}}$ denote the
set of $m\times n$ 0/1 tables with row
sums $\mathbf{r}$ and column sums $\mathbf{c}$.
Let $N=\sum_i r_i$ denote the number of edges in the corresponding bipartite graphs.

Our focus is on algorithms for estimating $|\Omega|$.
There are algorithms
\cite{JSV,BBV} for estimating $|\Omega|$ (and sampling (almost)
uniformly at random from $\Omega$)
which provably run in time polynomial in $n$ and~$m$ for {\it any\/} row/column sums.
We discuss these algorithms, which use Markov chain Monte Carlo (MCMC)
methods, in more detail later in the introduction.
In this paper, we study a simpler method known as
sequential importance sampling (SIS).

SIS has several purported advantages over
the more classical Markov chain Monte Carlo (MCMC) method, such as:
\begin{description}
\item[Speed:] Chen et al.~\cite{CDHL} claim that SIS is faster
than MCMC algorithms (their paper shows, by experiment, that for the
studied inputs, SIS is superior to the MCMC algorithm
of~\cite{BesagClifford}; moreover the authors state that they are
not aware of any MCMC-based algorithm that achieves similar results in
both accuracy and time as SIS).
In fact, Blanchet \cite{Blanchet} recently proved that SIS requires $O(N^2)$ time
when all of the row and column sums are at most $o(N^{1/4})$ (see
Bayati et al. \cite{BKS} for a related result for a different
algorithm). In contrast, we present a simple example where SIS
requires an exponentially large (in $n,m$) number of samples to
give an approximately correct answer. Note that, as mentioned
earlier, a MCMC algorithm was presented in \cite{JSV,BBV} which is
guaranteed to require at most time polynomial in $n,m$ for every
input. \item[Convergence Diagnostic:] One of the difficulties in
MCMC algorithms is determining when the Markov chain of interest
has reached the stationary distribution, unless we have analytical
bounds (as in the case of \cite{JSV,BBV}). SIS seemingly avoids
such complications since its output is guaranteed to be an
unbiased estimator of $|\Omega|$.  Unfortunately, it is unclear
how many estimates from SIS are needed before we have a guaranteed
close approximation of $|\Omega|$.  In our example for which SIS
requires exponential time, the estimator appears to converge, but
it converges to a quantity that is off from $|\Omega|$ by an
exponential factor.
\end{description}
Before formally stating our results, we detail the sequential importance sampling
approach for contingency tables,
following~\cite{CDHL}.
The general importance sampling paradigm
involves sampling from an `easy' distribution $\mu$ over
$\Omega$
that is, ideally, close to the uniform distribution.
At every round, the algorithm outputs a table $\mathbf{T}$ along with $\mu(\mathbf{T})$.
Since for any $\mu$ whose support is $\Omega$ one has
\[ E[1/\mu(\mathbf{T})] = |\Omega|,
\]
the algorithm takes many trials and outputs the average of
$1/\mu(\mathbf{T})$ as an estimate of $|\Omega|$.   More
precisely, let $\mathbf{T}^{(1)},\dots,\mathbf{T}^{(t)}$ denote
the outputs from $t$ trials of the SIS algorithm.  The
final estimate is
\begin{equation}
\label{estimator}
 X_t = \frac{1}{t}\sum_{\ell=1}^t \frac{1}{\mu(\mathbf{T}^{(\ell)})}.
\end{equation}
One typically uses a heuristic
to determine how many trials $t$ are needed until the
estimator has converged to the
desired quantity.

The sequential importance sampling algorithm of Chen et
al.~\cite{CDHL} constructs the table $\mathbf{T}$ in a
column-by-column manner. It is not clear how to order the columns
optimally, but this will not concern us as our negative results
will hold for any ordering of the columns. Suppose the procedure
is assigning column $i$ conditional on an existing assignment to columns $1,\dots,i-1$.
For $1\le j\le m$, let $r'_j$ be equal to
$r_j$ less the total number of $1$'s seen in row $j$ in columns
$1,\dots,i-1$.
Thus, $r'_1,\dots,r'_m$ are the residual row sums
after taking into account the assignments in the first $i-1$ columns.

The procedure of Chen et al.\ chooses column $i$ from the
following probability distribution. The
distribution is the projection onto column $i$
of the uniform distribution over assignments to
columns $i,\dots,n$ where the row sums are $r'_1,\dots,r'_m$
and column $i$ sums to $c_i$ (but ignoring the column sums
$c_{i+1},\dots,c_n)$.
The distribution is easy to describe
in closed form.
Let $\T_{1,i},\dots,\T_{m,i}\in \{0,1\}^m$ denote the assignment to column $i$,
where $\sum_j \T_{j,i}=c_i$.
Let $n'=n-i+1$ be the number of not yet assigned columns.
Clearly, $\T_{j,i}$ must be $0$ for every $j$ with $r'_j=0$, and $\T_{j,i}$ must be $1$ for every $j$ with $r'_j=n'$.
Let $J = \{j\in\{1,\dots,m\}~|~0<r'_j<n'\}$, that is, $J$ is the
set of rows whose entries are not forced to $0$ or $1$.
Then, the probability of the assignment $\T_{j,i}$ for $j\in J$ is proportional to

\begin{equation}
\label{eq:sisdistr}
  \prod_{j\in J}\left(\frac{r_j'}{n'-r_j'}\right)^{\T_{j,i}}.
\end{equation}
Sampling from this distribution over assignments
for column $i$ can be done efficiently by dynamic programming (see
Section 3.1 of~\cite{CDHL}).


\begin{note}
\label{ftn:delicatesampling} The described procedure may ``get
stuck'', that is, run into a situation when no valid assignment is
possible for the $i$-th column. In such case,
$1/\mu(\mathbf{T}^{(\ell)})$ is set to zero in \eqref{estimator}
for this trial and the procedure moves to the next trial. Chen et
al.\ also devised a more subtle sampling procedure for the $i$-th
column which never gets stuck. We do not describe this interesting
modification of the procedure, as the two procedures are
equivalent for the input instances which we discuss in this paper.
The reason is that for our instances even for the distribution given
by \eqref{eq:sisdistr} SIS never gets stuck.
\end{note}

We now state our negative result.  This is a simple family of examples
where the SIS
algorithm will grossly underestimate $|\Omega|$
unless the number of trials $t$ is exponentially large.
Our examples
will have the form $(1,1,\dots,1,d_r)$ for row sums and $(1,1,\dots,1,d_c)$
for column sums, where the number of rows is $m+1$, the number
of columns is $n+1$, and we require that $m+d_r=n+d_c$.
An important feature of our examples is that they are ``bad'' examples
regardless of whether the SIS procedure constructs the table column-by-column
or row-by-row.

\begin{theorem}
\label{thm:verybad}
Let $\beta>0,\gamma\in(0,1)$ be constants satisfying $\beta\neq\gamma$
and consider the input instances $\mathbf{r}=(1,1,\dots,1,\lfloor \beta m\rfloor)$,
$\mathbf{c}=(1,1,\dots,1,\lfloor \gamma m\rfloor)$ with $m+1$ rows.
Fix any order of columns (or rows, if sequential importance sampling
constructs tables row-by-row)
and let $X_t$ be the random variable
representing the estimate of the SIS procedure after $t$ trials of the algorithm,
that is, $X_t$ is given by ~\eqref{estimator}.
There exist constants $s_1\in(0,1)$ and $s_2>1$ such that
for every sufficiently large $m$ and for any $t\leq s_2^m$,
$$
\Prob{X_t\geq \frac{|\Omega_{\mathbf{r},\mathbf{c}}|}{s_2^m}  } \leq 3s_1^m.
$$
\end{theorem}

We mentioned earlier that there are MCMC algorithms which provably
run in time polynomial in $n$ and~$m$ for any row/column sums.
In particular, Jerrum, Sinclair, and Vigoda~\cite{JSV}
presented a polynomial time algorithm for estimating
the permanent of a non-negative matrix.  For the case of 0/1 matrices,
their result corresponds to a randomized algorithm, which for
a bipartite graph $G$,
estimates the number of perfect matchings of $G$ within a multiplicative
factor $(1\pm\eps)$ in time polynomial in $|G|$ and $1/\eps$.
The binary contingency tables problem studied in this paper
can be reduced to counting perfect matchings via a reduction of
Tutte \cite{Tutte}.
More recently, Bez\'akov\'a, Bhatnagar and Vigoda~\cite{BBV}
presented a related simulated annealing algorithm that works
directly with binary contingency tables to solve the problem
studied in this paper for all row/column sums, and has an improved
polynomial running time compared with \cite{JSV}.
We note that, in
addition to being formally asymptotically faster than any
exponential time algorithm, a polynomial time algorithm has
additional theoretical significance in that it (and its analysis)
implies non-trivial insight into the structure of the problem.

As a side note, we remark that even though SIS grossly underestimates the number of binary contingency tables for our examples with $m+1$ rows with row sums $(1,1,\dots,1,d_r)$ and $n+1=m+d_r-d_c+1$ columns with column sums $(1,1,\dots,1,d_c)$, it is possible to compute this number {\em exactly} using the formula $\binom{m}{d_c}\binom{n}{d_r}(m-d_c)! + \binom{m}{d_c-1}\binom{n}{d_r-1}(m-d_c+1)!$.

Some caveats are in order here.  Firstly, the above results imply only
that MCMC outperforms SIS asymptotically {\it in the worst case\/}; for
many inputs, SIS may well be much more efficient.  Secondly, the rigorous
worst case upper bounds on the running time of the above MCMC algorithms
are still far from practical.
Chen et al.~\cite{CDHL}
showed several examples where SIS outperforms MCMC methods.
We present a more systematic experimental study of the performance of SIS,
focusing on examples where all the row and column sums are identical as well
as on the ``bad'' examples from Theorem~\ref{thm:verybad}.
Our experiments suggest that SIS is extremely fast on the balanced examples,
while its performance on the bad examples confirms our theoretical analysis.
Understanding conditions under which SIS performs well  is, perhaps, the most
interesting open problem in the area. Specific problems include extending the result
of~\cite{Blanchet} to multiway contingency tables~\cite{CDS} and random graphs
with prescribed degrees~\cite{Blitzstein}.

We also note that the following simple modification of SIS may lead to
better performance.  Rather than assigning entries in a column-by-column
or row-by-row manner, assign at each step {\it either\/} the row
{\it or\/} the column with the largest residual sum.  It can easily be verified
that this enhanced scheme does produce correct results for the input
instances in Theorem~\ref{thm:verybad}.  However, we provide
experimental evidence that there are input instances for which even
this enhanced strategy fails.  These inputs are similar in flavor to those
in Theorem~\ref{thm:verybad}, but slightly more complicated.

We begin in Section \ref{sec:preliminaries}
by presenting a few basic lemmas that are used in the analysis of our negative example.
In Section \ref{sec:mainexample} we present  our main example where SIS is off by an
exponential factor, thus proving Theorem \ref{thm:verybad}.
Finally, in Section~\ref{sec:experiments} we summarize
some experimental results for SIS that support our theoretical analysis.

\section{Preliminaries}
\label{sec:preliminaries}
We will continue to let
 $\SISPr(\mathbf{T})$ denote the probability that a
table $\mathbf{T}\in\Omega_{\mathbf{r},\mathbf{c}}$ is generated by sequential importance
sampling algorithm.  We let $\TruePr(\mathbf{T})$ denote the uniform distribution
over $\Omega_{\mathbf{r},\mathbf{c}}$,
 which is the desired distribution.

Before beginning our main proofs we present two straightforward technical
lemmas which are used at the end of the proof of the main theorem.
The first lemma claims that if a large set of binary contingency
tables gets a very small probability under SIS, then
 SIS is likely to output an estimate which is not much bigger than
the size of the complement of this set, and hence very small.
For $S\subset\Omega_{\mathbf{r},\mathbf{c}}$,
let $\coS = \Omega_{\mathbf{r},\mathbf{c}}\setminus S$ denote its complement.

\begin{lemma}
\label{lem:underestimate} Let $p\leq 1/2$ and let
$S\subseteq\Omega_{r,c}$ be such that $\SISPr(S)\leq p$. Then for
any $a>1$, and any $t$, we have
\[    \Prob{ X_t\leq a |\coS|}  \geq 1 - pt
-2/a.
\]
\end{lemma}

\begin{proof}
The probability that all $t$ SIS trials are not in $S$ is at least
$$(1-p)^t  \geq 1-pt.$$
(This well-known inequality valid for $p\geq 0$ follows by induction on $t$.)

Let $\T^{(1)},\dots,\T^{(t)}$
 be the $t$ tables constructed by SIS. Then,
with probability at least $1-pt$,
we have $\T^{(\ell)}\in\coS$ for all $\ell$, $1\le
\ell\le t$.
Notice that for a table $\T$ constructed by SIS from $\coS$, we
have
\[
\ExpCond{\frac{1}{\SISPr(\T)}  }{ \T\in\coS}  =
\frac{|\coS|}{\mu(\coS)}.
\]
Let $\F$ denote the event that $\T^{(\ell)}\in\coS$ for all $\ell$, $1\le
\ell\le t$; hence,
\[  \ExpCond{X_t}{\F} =  \frac{|\coS|}{\mu(\coS)}.
\]

We can use Markov's inequality to estimate the probability that
SIS returns an answer which is more than a factor of $a$ worse
than the expected value, conditioned on the fact that no SIS trial
is from $S$:
$$
\ProbCond{X_t>a|\coS| }{\F} \leq
\ProbCond{X_t>(a/2)\frac{|\coS|}{\mu(\coS)} }{\F} \leq
\frac{2}{a},$$ where in the first inequality we used
$\mu(\coS)\geq 1/2$.

Finally, removing the conditioning we get:
\begin{eqnarray*}
\Prob{X_t\leq a |\coS|} &\geq& \ProbCond{ X_t\leq a
|\coS|}{\F}\Prob{\F}
\\
& \geq & \left(1-\frac{2}{a}\right)(1-pt)
\\
&\geq & 1-pt-\frac{2}{a}.
\end{eqnarray*}
\end{proof}

The second technical lemma shows that if in a row with large sum
(linear in $m$) there exists a large number of columns (again linear in $m$)
for which the SIS probability of
placing a $1$ at the corresponding position differs significantly
from the correct probability,
then in any subexponential number of trials
the SIS estimator will very likely exponentially underestimate the correct answer.

Let $\mathcal{A}_{i-1}$ denote the set of all assignments of
$0/1$ to columns $1,\dots,i-1$ such that the column sums are satisfied and
none of the row sums are exceeded.
Thus, $\mathbf{A}\in\mathcal{A}_{i-1}$
denotes that $\mathbf{A}$ is a specific assignment of $0/1$ to
the first $i-1$ columns.    Finally, for $\mathbf{A}\in\mathcal{A}_{i-1}$,
we use the following notation:
\[ \TruePr\left(\T_{j,i} = 1 \, \vert \, \T_{i-1} = \mathbf{A}\right) =
\frac{\TruePr\left(S''\right)}{\TruePr\left(S'\right)},
\]
where \[ S' = \{\T\in\Orc: \mbox{the first $i-1$ columns of $\T$
are the same as $\mathbf{A}$}\},
\]
and
\[ S'' = \{\T\in S': \T_{j,i} = 1 \}.
\]
Similarly, we use
\begin{equation}\label{krtko2}
\SISPr\left(\T_{j,i} = 1 \, \vert \, \T_{i-1} = \mathbf{A}\right) = \frac{\mu(S'')}{\mu(S')}.
\end{equation}

As we mentioned in Remark \ref{ftn:delicatesampling}, SIS will never get stuck for our
input instances. For such input instances, \eqref{krtko2} is the
same as the probability that SIS assigns 1 to $\T_{j,i}$,
given
that the first $i-1$ columns are filled with $\mathbf{A}$.

\begin{lemma}
\label{lem:separable}
Let $\alpha<\beta$ be positive constants.
Consider a class of instances of the binary contingency tables
problem, parameterized by $m$,
with $m+1$ row sums, the last of which is $\lfloor\beta m\rfloor$.
The remaining row sums and column sums can be arbitrary as long as
the SIS procedure never gets stuck.
Suppose that there exist constants
$f<g$ and a set $I$ of cardinality $\lfloor\alpha m\rfloor >0$
such that one of the following statements is true:
\begin{enumerate}[(i)]
\item
for every $i\in I$ and any $\mathbf{A}\in\mathcal{A}_{i-1}$,
$$
\TruePr(\T_{m+1,i}=1~|~\T_{i-1}=\mathbf{A})\leq f < g \leq \SISPr(\T_{m+1,i}=1~|~\T_{i-1}=\mathbf{A}),
$$
\item
for every $i\in I$ and any $\mathbf{A}\in\mathcal{A}_{i-1}$,
$$
\SISPr(\T_{m+1,i}=1~|~\T_{i-1}=\mathbf{A})\leq f < g \leq \TruePr(\T_{m+1,i}=1~|~\T_{i-1}=\mathbf{A}).
$$
\end{enumerate}
Then there exists a constant $b_1\in (0,1)$ such that for any constant $1<b_2<1/b_1$
and any sufficiently large $m$,
for any $t\leq b_2^m$,
$$\Prob{X_t \geq \frac{|\Omega_{\mathbf{r},\mathbf{c}}|}{b_2^{m} }  }   \leq 3(b_1b_2)^m.$$
\end{lemma}


\begin{proof}
We will analyze case (i); the other case follows from
analogous arguments.

Let $U_1,\dots,U_n$ be the entries in the last row of a uniformly
random contingency table with the prescribed row and column sums.
(Note that $U_1,\dots,U_n$ are random variables.) Similarly, let
$V_1,\dots,V_n$ be the entries in the last row of a contingency
table with the prescribed row and column sums generated by SIS.

The random variable $\X_i$ is dependent on $\X_j$ for $j<i$ and
$\Y_i$ is dependent on $\Y_j$ for $j<i$. However, for every $i\in I$, $\X_i$ is stochastically dominated by $\X'_i$, where
$\X'_i$, $i\in I$ is a set of independent Bernoulli random variables that take value $1$ with probability $f$.
Similarly, for every $i\in I$, $\Y_i$ stochastically dominates $\Y'_i$, where
$\Y'_i$, $i\in I$ is a set of independent Bernoulli random variables that take value $1$ with probability $g$.


Now we may use the Chernoff bound
(see, e.\.g.,~\cite{DP09}, Theorem 1.1).
 Let $k=\lfloor \alpha
m\rfloor$. Then
  $$\Prob{   \sum_{i\in I}\X_i' - kf > \frac{g-f}{2}k  }
  \leq \exp({-(g-f)^2k/2})$$
  and
  $$
  \Prob{kg-\sum_{i\in I}\Y_i' > \frac{g-f}{2}k }
  \leq \exp({-(g-f)^2k/2}).
  $$
Let $S$ be the set of all tables which have less than
$kf+(g-f)k/2 = kg-(g-f)k/2$
ones in the last row of the columns in $I$. Let
$b_1:=\exp({-(g-f)^2\alpha/4})\in(0,1)$. Then $\exp({-(g-f)^2k/2})\leq b_1^m$
for $m\geq 1/\alpha$. Thus, by the first
inequality, under the uniform distribution over all binary contingency tables the
probability of the set $S$ is at least
$1-b_1^m$. However, by the second
inequality, SIS constructs a table
from the set $S$ with probability at most
$b_1^m$.

We are ready to use Lemma~\ref{lem:underestimate} with $S$ as
defined above and
$p=b_1^m$. Since under the uniform distribution the probability of $S$
is at least $1-b_1^m$, we have that $|S|\geq
(1-b_1^m)|\Omega_{\mathbf{r},\mathbf{c}}|$. Let $b_2\in(1,1/b_1)$
be any constant and consider $t\leq b_2^m$ SIS trials. Let
$a=(b_1b_2)^{-m}$. Then, by Lemma~\ref{lem:underestimate}, with
probability at least $1-pt-2/a\geq 1-3(b_1b_2)^m$ the SIS
procedure outputs a value which is at most an $ab_1^m=b_2^{-m}$
fraction of $|\Omega_{\mathbf{r},\mathbf{c}}|$.
\end{proof}

\section{Proof of Main Theorem}

\label{sec:mainexample}

In this section we prove Theorem \ref{thm:verybad}.
Before we analyze the input instances from Theorem \ref{thm:verybad},
we first consider
the following simpler class of inputs.

\subsection{Row sums $(1,1,\dots,1,d)$ and column sums $(1,1,\dots,1)$}

The row sums are
$(1,\dots,1,d)$ and the number of rows is $m+1$.
The column sums are $(1,\dots,1)$ and the number of columns is
$n=m+d$.
We assume that sequential importance sampling
constructs the tables column-by-column.
If SIS constructed the tables row-by-row, starting
with the row with sum $d$, then it would in fact output the
correct number of tables exactly.
However, in the next subsection
we will use this simplified case as a tool in our analysis of the
input instances $(1,\dots,1,d_r)$, $(1,\dots,1,d_c)$, for which SIS must necessarily
fail regardless of whether it works row-by-row or
column-by-column, and regardless of the order it chooses.

\begin{lemma}
\label{lem:columnones} Let $\beta>0$, and consider an input of the form
$(1,\dots,1,d),(1,\dots,1)$
with $m+1$ rows where $d=\lfloor\beta m\rfloor$.
Then there exists a constants $s>1$, such that for any sufficiently
large $m$,
for any $t\leq s^m$,
$$\Prob{X_t\geq \frac{|\Omega_{\mathbf{r},\mathbf{c}}|}{s^{m}}  }   \leq 3(1/2)^m.$$
\end{lemma}

The idea for the proof of the lemma is straightforward. By the
symmetry of the column sums, for large $m$ and $d$ and
$\alpha\in(0,1)$ a uniform random table will have about
$\alpha d$ ones in the first $\alpha n$ cells of the last row,
with high probability. We will show that for some $\alpha\in(0,1)$
and $d=\beta m$,  sequential importance sampling is very
unlikely to put this many ones in the first $\alpha n$
columns of the last row. Therefore, since with high probability
sequential importance sampling will not construct any table
from a set that is a large fraction of all legal tables, it will likely
drastically underestimate the number of tables.

Before we prove the lemma, let us first
compare the column distributions arising from the uniform distribution
over all binary contingency tables with the SIS
distributions. We refer to the column distributions induced by
the uniform distribution over all tables as the {\em true}
distributions.
The true probability of $1$ in the first column and
last row can be computed as the number of tables with $1$ at this
position divided by the total number of tables. For the
sequence of row and column sums specified in the
statement of Lemma \ref{lem:columnones},
let $Z(m,d)$ denote the total number of tables with these row/column sums.
Note, $Z(m,d)={n \choose d}m! = {m+d \choose d}m!$, since a table is uniquely specified by the positions of
ones in the last row and the permutation matrix in the remaining
rows and corresponding columns. Therefore,
$$\TruePr(\T_{m+1,1}=1) = \frac{Z(m,d-1)}{Z(m,d)} =\frac{{m+d-1 \choose d-1}m!}{{m+d\choose
d}m!} = \frac{d}{m+d}.$$

On the other hand, by the definition of sequential importance
sampling, $\Prob{\T_{i,1}=1}\propto r_i/(n-r_i)$, where $r_i$ is the row
sum in the $i$-th row. Therefore,
$$\SISPr(\T_{m+1,1}=1) = \frac{\frac{d}{n-d}}{\frac{d}{n-d} + m
\frac{1}{n-1}} = \frac{d(m+d-1)}{d(m+d-1)+m^2}. $$

Observe that if $d \approx \beta m$ for some constant $\beta>0$,
then for sufficiently large $m$ we have
\[
\SISPr(\T_{m+1,1}=1) > \TruePr(\T_{m+1,1}=1).
\] As we will see, this will
be true for a linear number of columns, which turns out to be
enough to prove that in polynomial time sequential importance sampling
exponentially underestimates the total number of
binary contingency tables with high probability.

\begin{proof}[Proof of Lemma~\ref{lem:columnones}]
We will find a constant $\alpha$ such that for every column
$i<\alpha m$ we will be able to derive an upper bound on the true
probability and a lower bound on the SIS probability of
$1$ appearing at the $(m+1,i)$ position.

For a partially filled table with columns $1,\dots,i-1$ assigned,
let $d_i$ be the remaining sum in the last row and let $m_i$ be the
number of other rows with remaining row sum $1$ (note that this determines
the contents $\mathbf{A}\in\mathcal{A}_{i-1}$ of the first $i-1$ columns, up to permutation).
Then the true probability of $1$ in the $i$-th column and last row can be bounded as
  $$\TruePr(\mathbf{T}_{m+1,i}=1~|~\mathbf{T}_{i-1}=\mathbf{A})=\frac{d_i}{m_i+d_i}
  \leq \frac{d}{m+d-(i-1)} =: f(d,m,i),
  $$
while the probability under SIS can be bounded as
  \begin{eqnarray*}
  \SISPr(\mathbf{T}_{m+1,i}=
  1~|~\mathbf{T}_{i-1}=\mathbf{A})
  &=&
  \frac{d_i(m_i+d_i-1)}{d_i(m_i+d_i-1)+m_i^2}
  \\
  &\geq&
  \frac{(d-(i-1))(m+d-i)}{d(m+d-1)+m^2}
  \\ & =: & g(d,m,i).
  \end{eqnarray*}
Observe that for fixed $m,d$, the function $f$ is increasing and the function
$g$ is decreasing in $i$, for $i<d$.

Recall that we are considering a family of input instances
parameterized by $m$ with $d=\lfloor\beta m\rfloor$,
for a fixed $\beta>0$. We will
consider $i<\alpha m$ for some $\alpha\in(0,\beta)$. Let
\begin{equation}
\label{eq:finf}
f^{\infty}(\alpha,\beta):=\lim_{m\to\infty} f(d,m,\alpha m) =
\frac{\beta}{1+\beta-\alpha};
\end{equation}
\begin{equation}
\label{eq:ginf}
g^{\infty}(\alpha,\beta):=\lim_{m\to\infty} g(d,m,\alpha m) =
\frac{(\beta-\alpha)(1+\beta-\alpha)}{\beta(1+\beta)+1};
\end{equation}
\begin{equation}
\triangle_{\beta}:=g^{\infty}(0,\beta)-f^{\infty}(0,\beta)=\frac{\beta^2}
{(1+\beta)(\beta(1+\beta)+1)}>0,
\end{equation}
and observe that for fixed $\beta$,
$f^{\infty}$ is increasing in $\alpha$
and $g^{\infty}$ is decreasing in $\alpha$, for $\alpha<\beta$.
Let $\alpha$, $0<\alpha < \beta$ be such that $g^{\infty}(\alpha,\beta)-f^{\infty}(\alpha,\beta)\geq
\triangle_{\beta}/2$.
Such an $\alpha$ exists by continuity (we only need to take a small enough $\alpha$).

By the above, for any $\eps>0$ and sufficiently large $m$, and for
any $i<\alpha m$, the true probability is
upper-bounded by $f^{\infty}(\alpha,\beta)+\eps$
and the SIS probability is lower-bounded by $g^{\infty}(\alpha,\beta)-\eps$.
For our purposes it is enough to fix $\eps=\triangle_{\beta}/8$.
Now we can use Lemma~\ref{lem:separable} with $\alpha$ and
$\beta$ defined as above, $f=f^{\infty}(\alpha,\beta)+\eps$ and
$g=g^{\infty}(\alpha,\beta)-\eps$ (notice that all these constants
depend only on $\beta$), and $I=\{1,\dots,\lfloor\alpha m\rfloor\}$.
Let $b_1\in (0,1)$ be the constant guaranteed by Lemma~\ref{lem:separable}
and let $b_2=1/(2b_1)$. This finishes the proof of the lemma with $s=b_2$.
%
%
\end{proof}

\begin{note}
Notice that every contingency table with row sums $(1,1,\dots,1,d)$
and column sums $(1,1,\dots,1)$ is binary. Thus, this instance
proves that the column-based SIS procedure for general
(non-binary) contingency tables~\cite{CDHL} has the same flaw as the binary
SIS procedure.
We expect
that the negative
example used for Theorem~\ref{thm:verybad} also extends to general (\ie, non-binary)
contingency tables,
but the analysis becomes more cumbersome.
\end{note}

\subsection{Row sums $(1,1,\dots,1,d_r)$ and column sums $(1,1,\dots,d_c)$}



We will now prove our main result, using ideas from the proof of Lemma~\ref{lem:columnones}.

\begin{proof}[Proof of Theorem~\ref{thm:verybad}]
Recall that we are working with row sums $(1,1,\dots,1,d_r)$,
where the number of rows is $m+1$, and column sums
$(1,1,\dots,1,d_c)$,
where the number of columns is $n+1=m+1+d_r-d_c$. We will
eventually fix $d_r=\lfloor \beta m\rfloor$ and $d_c=\lfloor \gamma
m\rfloor$, but to simplify our expressions we work
with $d_r$ and $d_c$ for now.

The
theorem claims that the SIS procedure fails for
an arbitrary order of columns with
high probability. We first analyze the case when the SIS procedure starts with columns
of sum $1$; we shall address the issue of arbitrary column order later.
As before, under the assumption that the first column has sum $1$,
we compute
the probabilities of 1 being in the last row for uniform random tables
and for SIS respectively.
For the true probability, the total number of tables can be
computed as ${m \choose d_c}{n \choose
d_r}(m-d_c)! + {m \choose d_c-1}{n \choose
d_r-1}(m-d_c+1)!$, since a table is uniquely determined by the
positions of ones in the $d_c$ column and $d_r$ row and a
permutation matrix on the remaining rows and columns.
Thus we have
  \begin{eqnarray*}
  \TruePr(\T_{m+1,1}=1)
  &=&
  \frac{{m\choose d_c}{n-1\choose d_r-1}(m-d_c)! +
  {m\choose d_c-1}{n-1\choose d_r-2}(m-d_c+1)!}{{m \choose d_c}{n \choose
  d_r}(m-d_c)! + {m \choose d_c-1}{n \choose
  d_r-1}(m-d_c+1)!}
  \\
  &=&
  \frac{d_r(n-d_r+1) + d_cd_r(d_r-1)}{n(n-d_r+1)+nd_cd_r}
  =: f_2(m,d_r,d_c);\\
  \\
  \SISPr(\T_{m+1,1}=1) &=& \frac{\frac{d_r}{n-d_r}}{\frac{d_r}{n-d_r}+m\frac{1}{n-1}}
  = \frac{d_r(n-1)}{d_r(n-1)+m(n-d_r)} =: g_2(m,d_r,d_c).
  \end{eqnarray*}
 Let $d_r = \lfloor\beta m\rfloor$ and $d_c = \lfloor\gamma m\rfloor$
for some constants $\beta>0,\gamma\in (0,1)$ (notice that this
choice guarantees that $n\geq d_r$ and $m\geq d_c$, as required).
Then, as $m$ tends to infinity, $f_2$ approaches
  $$f^{\infty}_2(\beta,\gamma):=\frac{\beta}{1+\beta-\gamma},$$
and $g_2$ approaches
  $$g^{\infty}_2(\beta,\gamma):=\frac{\beta(1+\beta-\gamma)}
  {\beta(1+\beta-\gamma)+1-\gamma}.$$
Notice that $f^{\infty}_2(\beta,\gamma)=
g^{\infty}_2(\beta,\gamma)$ if and only if $\beta=\gamma$.
Moreover, $f_2^\infty(\beta, \gamma) < g_2^\infty(\beta, \gamma)$
if and only if $\beta > \gamma$.
Suppose that $\beta>\gamma$, that is,
$f^{\infty}_2(\beta,\gamma)<g^{\infty}_2(\beta,\gamma)$ (the opposite
case follows analogous arguments and uses the second part
of Lemma~\ref{lem:separable}).
As in the proof of Lemma~\ref{lem:columnones}, we can
define $\alpha$ such that if the importance sampling does not
choose the column with sum $d_c$ in its first $\alpha m$ choices,
then in any subexponential number of trials it will exponentially
underestimate the total
number of tables with high
probability. Formally, we derive an upper bound on the true
probability of $1$ being in the last row of the $i$-th column, and
a lower bound on the SIS probability of the same event (both
conditioned on the fact that the $d_c$ column is not among the
first $i-1$ columns assigned). Assume that we already assigned
the first $i-1$ columns of the table. Let $d_r^{(i)}$ be the current
residual sum in the last row (that is, $d_r^{(i)}$ is $d_r$ less the
number of ones assigned to the last row of columns $1,\dots,i-1$),
$m_i$ be the remaining number of rows with sum $1$, and $n_i$ the remaining number of
columns with sum $1$ (note that this determines the contents $\mathbf{A}\in\mathcal{A}_{i-1}$ of the first
$i-1$ columns, up to permutation). Notice that
$n_i=n-i+1$, $m\geq m_i\geq m-i+1$, and $d_r\geq d_r^{(i)}\geq
d_r-i+1$. Then
  \begin{eqnarray*}
  \TruePr(\mathbf{T}_{m+1,i}=1~|~\mathbf{T}_{i-1}=\mathbf{A})
  &=&
  \frac{d_r^{(i)}(n_i-d_r^{(i)}+1) + d_cd_r^{(i)}(d_r^{(i)}-1)}
  {n_i(n_i-d_r^{(i)}+1)+n_id_cd_r^{(i)}}
  \\
  &\leq&
  \frac{d_r(n-d_r+1) + d_cd_r^2}
  {(n-i+1)(n-i-d_r+2) + (n-i+1)d_c(d_r-i+1)}
  \\
  &\leq&
  \frac{d_r(n-d_r+1) + d_cd_r^2}
  {(n-i)(n-i-d_r) + (n-i)d_c(d_r-i)}
 \\
  &=: & f_3(m,d_r,d_c,i);
  \\
  \SISPr(\mathbf{T}_{m+1,i}=1~|~\mathbf{T}_{i-1}=\mathbf{A})
  &=& \frac{d_r^{(i)}(n_i-1)}{d_r^{(i)}(n_i-1)+m_i(n_i-d_r^{(i)})}\\
  &\geq&  \frac{(d_r-i)(n-i)}{d_rn+m(n-d_r)}
  \\
&  =:& g_3(m,d_r,d_c,i).
  \end{eqnarray*}
As before, notice that if we fix $m, d_r, d_c>0$ satisfying
$d_c<m$ and $d_r<n$, then $f_3$ is
an increasing function and $g_3$ is a decreasing function in
$i$, for $i<\min\{n-d_r,d_r\}$. Recall that $n-d_r=m-d_c$.

Let $\alpha$ be a number such that $0<\alpha<\min\{1-\gamma,\beta\}$
(we will further specify how $\alpha$ is chosen shortly---it will be small enough
to satisfy equations~\eqref{eq:krt3} and~\eqref{eq:krt4} below).
Suppose that $i\leq \alpha m<\min\{m-d_c,d_r\}$.
Thus, the upper bound on $f_3$ in this range of $i$ is $f_3(m,d_r,d_c,\alpha m)$
and the lower bound on $g_3$ is $g_3(m,d_r,d_c,\alpha m)$.
If $d_r=\lfloor\beta m\rfloor$ and $d_c=\lfloor\gamma m\rfloor$,
then the upper bound on $f_3$ converges to
\begin{equation}\label{eq:krt1}
f^{\infty}_3(\alpha,\beta,\gamma) := \lim_{m\to\infty}f_3(m,d_r,d_c,\alpha m) =
\frac{\beta^2}
{(1+\beta-\gamma-\alpha)(\beta-\alpha)}
\end{equation}
and the lower bound on $g_3$ converges to
\begin{equation}\label{eq:krt2}
g^{\infty}_3(\alpha,\beta,\gamma) := \lim_{m\to\infty}g_3(m,d_r,d_c,\alpha m) =
\frac{(\beta-\alpha)(1+\beta-\gamma-\alpha)}{\beta(1+\beta-\gamma)+1-\gamma}.
\end{equation}
Let
$$\triangle_{\beta,\gamma}:=g^{\infty}_3(0,\beta,\gamma)-f^{\infty}_3(0,\beta,\gamma)
=g^{\infty}_2(\beta,\gamma)-f^{\infty}_2(\beta,\gamma)>0.$$
We set $\alpha>0$ to satisfy
\begin{equation}\label{eq:krt3}
g^{\infty}_3(\alpha,\beta,\gamma)-f^{\infty}_3(\alpha,\beta,\gamma)\geq
\triangle_{\beta,\gamma}/2.
\end{equation}
(Small enough $\alpha$ will work, by the continuity of \eqref{eq:krt1} and
\eqref{eq:krt2} for $\alpha\in(0,\min\{1-\gamma,\beta\})$.

Now we can conclude this part of the
proof identically to the last paragraph of the proof of
Lemma~\ref{lem:columnones}.

It remains to deal with the case when sequential importance
sampling picks the $d_c$ column within the first $\lfloor\alpha m\rfloor$
columns. Suppose $d_c$ appears as the $k$-th column.
In this case we focus on the subtable consisting of the last
$n+1-k$ columns with sum $1$,
$m'$ rows with sum $1$, and one row with sum~$d'$,
an instance of the form
$(1,1,\dots,1,d'),(1,\dots,1)$. We will use arguments similar to the proof of
Lemma~\ref{lem:columnones}.

First we express $d'$ as a function of $m'$.
The number of rows with row sum $1$ decreased by at least $d_c-1=\lfloor \gamma m\rfloor-1\geq \gamma m-2$,
and at most by $(\alpha+\gamma)m$.
Hence, $(1-\alpha-\gamma)m\leq m'\leq (1-\gamma)m+2$. Similarly, $d_r-\alpha m\leq d'\leq d_r$ where
$d_r=\lfloor\beta m\rfloor\geq\beta m -1$. Let $\beta'$ be such that $d'=\beta'm'$.
Thus, $(\beta-\alpha-1/m)/(1-\gamma+2/m)\leq\beta'\leq \beta/(1-\alpha-\gamma)$.

Now we find $\alpha'$ such that for any $i\leq \alpha'm'$ we will
be able to derive an upper bound on the true probability and a
lower bound on the SIS probability of $1$ appearing at position $(m'+1,i)$
of the $(n+1-k)\times m'$ subtable, no matter how the first $k$ columns were
assigned.

By the derivation in the proof
of Lemma~\ref{lem:columnones} (see expressions
\eqref{eq:finf} and \eqref{eq:ginf}),
as $m'$ (and thus also $m$) tends to infinity,
the upper bound on the true probability approaches
\begin{eqnarray*}
f^{\infty}(\alpha',\beta')&=&\lim_{m\to\infty}\frac{\beta'}{1+\beta'-\alpha'}
\\
&\leq &\lim_{m\to\infty}
\frac{\frac{\beta}{1-\alpha-\gamma}}{1+\frac{\beta-\alpha-\frac{1}{m}}{1-\gamma+\frac{2}{m}}-\alpha'}
\\
&=& \frac{\frac{\beta}{1-\alpha-\gamma}}{1+\frac{\beta-\alpha}{1-\gamma}-\alpha'}
\\
&=:&
f^{\infty}_4(\alpha,\beta,\gamma,\alpha')
\end{eqnarray*}
and the lower bound on the SIS probability approaches
\begin{eqnarray*}
g^{\infty}(\alpha',\beta')&=&\lim_{m\to\infty}\frac{(\beta'-\alpha')(1+\beta'-\alpha')}{\beta'(1+\beta')+1}
\\ &\geq&
\lim_{m\to\infty}\frac{(\frac{\beta-\alpha-\frac{1}{m}}{1-\gamma+\frac{2}{m}}-\alpha')(1+\frac{\beta-\alpha-\frac{1}{m}}{1-\gamma+\frac{2}{m}}-\alpha')}{\frac{\beta}{1-\alpha-\gamma}(1+\frac{\beta}{1-\alpha-\gamma})+1}\\
&=&\frac{(\frac{\beta-\alpha}{1-\gamma}-\alpha')(1+\frac{\beta-\alpha}{1-\gamma}-\alpha')}{\frac{\beta}{1-\alpha-\gamma}(1+\frac{\beta}{1-\alpha-\gamma})+1}
\\
&=: &
g^{\infty}_4(\alpha,\beta,\gamma,\alpha').
\end{eqnarray*}

Let us evaluate $f^{\infty}_4$ and $g^{\infty}_4$ for $\alpha=\alpha'=0$:
$$f^{\infty}_4(0,\beta,\gamma,0) = \frac{\frac{\beta}{1-\gamma}}{1+\frac{\beta}{1-\gamma}}$$
and
$$g^{\infty}_4(0,\beta,\gamma,0) = \frac{\frac{\beta}{1-\gamma}(1+\frac{\beta}{1-\gamma})}{\frac{\beta}{1-\gamma}(1+\frac{\beta}{1-\gamma})+1}.$$
Substituting $x$ for $\beta/(1-\gamma)$, we can see that $f^{\infty}_4(0,\beta,\gamma,0) < g^{\infty}_4(0,\beta,\gamma,0)$ since $\beta/(1-\gamma)\geq 0$.

Now let
$\triangle'_{\beta,\gamma}:=g^{\infty}_4(0,\beta,\gamma,0)-f^{\infty}_4(0,\beta,\gamma,0)>0$. By continuity,
for small enough $\alpha,\alpha'>0$ we have
\begin{equation}\label{eq:krt4}
g^{\infty}_4(\alpha,\beta,\gamma,\alpha')-f^{\infty}_4(\alpha,\beta,\gamma,\alpha')\geq\triangle'_{\beta,\gamma}/2.
\end{equation}

Now we proceed in a fashion similar to the last paragraph of the proof of
Lemma~\ref{lem:columnones}. More precisely, let $\eps:=\triangle'_{\beta,\gamma}/8$ and
let $f:=f^{\infty}_4(\alpha,\beta,\gamma,\alpha')+\eps$ and
$g:=g^{\infty}_4(\alpha,\beta,\gamma,\alpha')-\eps$ be the upper bound
(for sufficiently large $m$) on
the true probability and the lower bound
on the SIS probability of
$1$ appearing at the position $(m+1,i)$ for $i\in
I:=\{k+1,\dots,k+\lfloor\alpha'm'\rfloor\}$.
Therefore Lemma~\ref{lem:separable} with parameters $\alpha'$,
$\beta$, $I$ of size $|I|=\lfloor\alpha'm'\rfloor$,
$f$, and $g$ implies the statement of the theorem.

Finally, if the SIS procedure constructs the tables row-by-row instead of
column-by-column, symmetrical arguments hold. This completes the
proof of Theorem~\ref{thm:verybad}.
\end{proof}

\section{Experiments}
\label{sec:experiments}

We performed several experimental tests which show sequential
importance sampling to be a promising approach for certain classes
of input instances.

We ran the sequential importance sampling algorithm for binary contingency
tables, using the following stopping heuristic. Let $N=n+m$.
For some $\eps,k>0$ we stopped if the last $kN$ estimates were all
within a $(1+\eps)$ factor of the current estimate. We set $\eps=0.01$
and $k=5$.

\begin{figure}[H]
  \begin{minipage}[t]{.5\textwidth}
    \begin{center}
      \includegraphics[type=eps,ext=.eps,read=.eps,scale=0.8]{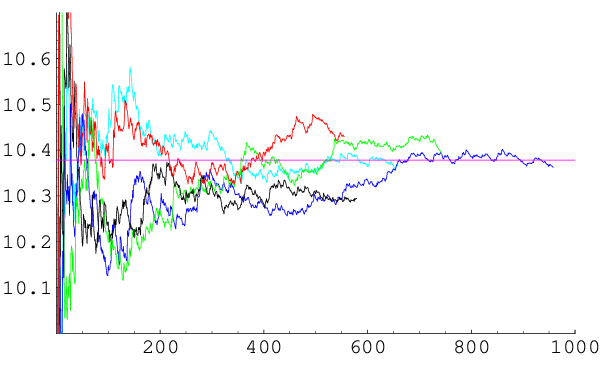}
      (a)
    \end{center}
  \end{minipage}
  \hfill
  \begin{minipage}[t]{.5\textwidth}
    \begin{center}
      \includegraphics[type=eps,ext=.eps,read=.eps,scale=0.8]{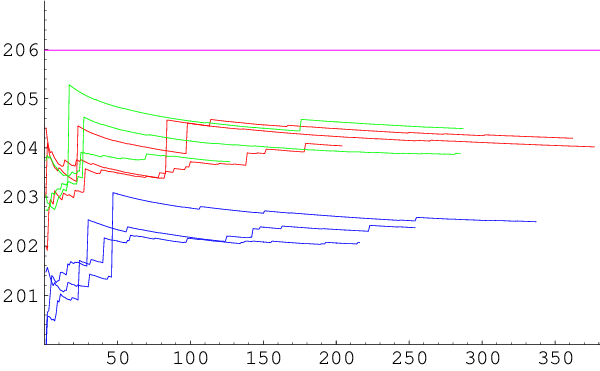}
      (b)
    \end{center}
  \end{minipage}
      \caption{
      The estimate produced by sequential importance sampling as a function
      of the number of trials on two different instances. In both figures, the
      horizontal line shows the correct number of
      corresponding binary contingency tables.
      (a) The left instance is a $50\times 50$ matrix where all $r_i=c_j=5$.
      The $x$-axis is the number of SIS trials, and the
       $y$-axis corresponds to the estimate scaled down by a factor of
      $10^{280}$. Five independent runs
      of sequential importance sampling are depicted.
      Notice that the $y$-axis ranges from $10$ to $10.7$, a relatively
      small interval, thus it appears SIS converges to the correct estimate.
      (b) The input instance is from Theorem~\ref{thm:verybad}
      with $m=300$, $\beta=0.6$ and $\gamma=0.7$.
      The estimate ($y$-axis) is plotted on a logarithmic scale (base $10$)
      and one unit on the $x$-axis corresponds to $1000$ SIS trials.
      Note that in this
      instance SIS appears to converge to an incorrect estimate.
      Nine independent runs of the SIS algorithm are shown:
      the red curves construct tables column-by-column with columns
      sorted by decreasing sum, the blue curves construct row-by-row with
      rows sorted by decreasing sum, and the green curves construct column-by-column
      with columns sorted increasingly.}
      \label{fig:convergence}
\end{figure}

Figure \ref{fig:convergence}(a) shows the evolution of the SIS
estimate as a function of the number of trials on the input with all row and column sums
$r_i=c_j=5$, and $50\times 50$ matrices. In our simulations we
used the more delicate sampling mentioned in Remark
\ref{ftn:delicatesampling}, which guarantees that the assignment in every
column is valid, \ie, such an assignment can always be extended to a valid table
(or, equivalently, that the random variable $X_t$ is always strictly positive).
Five independent
runs are depicted, together with the correct number of tables
$ \approx 1.038\times 10^{281}$, which we computed exactly. To make the figure legible,
the $y$-axis is scaled by a factor of $10^{280}$ and it only shows
the range from $10$ to $10.7$.
Note that the algorithm appears to converge to the correct estimate, and our
stopping heuristic appears to capture this behavior.

In contrast, Figure \ref{fig:convergence}(b) depicts the SIS evolution
on the negative example from Theorem \ref{thm:verybad}
with $m=300, \beta=0.6$ and $\gamma=0.8$, \ie, the
input is $(1,\dots,1,179),(1,\dots,1,240)$ on a $301\times 240$
matrix.
In this case the correct number of tables is
\[
{300\choose 240}{239\choose 179}(300-240)! + {300\choose 239}{239 \choose 178}(300-239)!
  \approx 9.684\times 10^{205}.
  \]
We ran the SIS algorithm under three different settings:
first, we constructed the tables column-by-column where the
columns were ordered from the largest sum, as suggested in the
paper by
Chen et al.~\cite{CDHL}
(the red curves correspond to three
independent runs with this setting); second,
we ordered the columns from the smallest sum (the green curves);
and third, we constructed the tables row-by-row where
the rows were ordered from the largest sum
(the blue curves).
The $y$-axis is on a {\em logarithmic scale} (base 10) and one unit on the
$x$-axis corresponds to $1000$ SIS trials. We ran the SIS estimates
for twice the number of trials determined by our stopping heuristic to
indicate that the
unfavorable performance of the SIS estimator on this
example is not the result of a poor choice of stopping heuristic.
Notice that even the best estimator differs from the true value by
about a factor of $40$, while the blue curves are off by more than
a factor of $1000$.

\begin{figure}[H]
\begin{center}
  \includegraphics[type=eps,ext=.eps,read=.eps,scale=1]{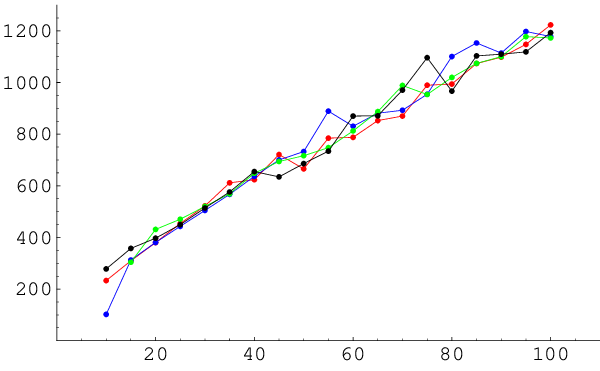}
\end{center}
      \caption{The number of SIS trials before the algorithm
      converges,
      as a function of the input size.
      The curves correspond to $5$ (red), $10$ (blue), $\lfloor 5\log n\rfloor$ (green), and $\lfloor
n/2\rfloor$ (black) regular row and column sums.}
      \label{fig:compareReg}
\end{figure}

\begin{figure}[H]
  \begin{minipage}[t]{.45\textwidth}
    \begin{center}
      \includegraphics[type=eps,ext=.eps,read=.eps,scale=.68]{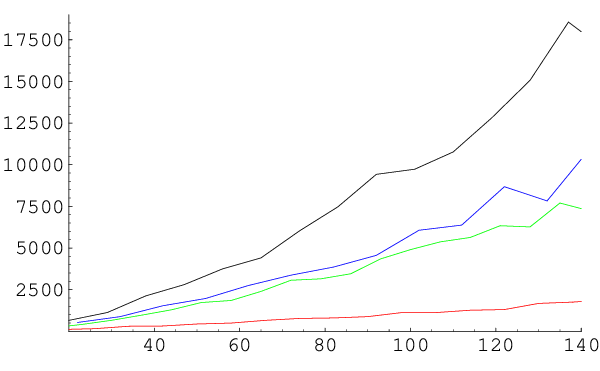}
    \end{center}
  \end{minipage}
  \hfill
  \begin{minipage}[t]{.45\textwidth}
    \begin{center}
      \includegraphics[type=eps,ext=.eps,read=.eps,scale=.65]{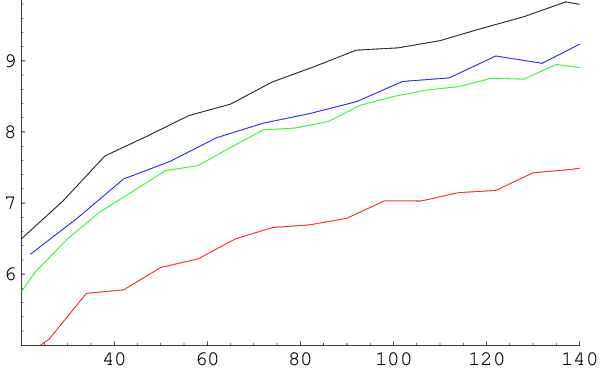}
    \end{center}
  \end{minipage}

\caption{The number of SIS trials until the algorithm converges as
      a function of $m+n$.  The inputs are
      of the type described in Theorem~\ref{thm:verybad}, with
      $\beta=0.1, \gamma=0.5$ (red), $\beta=\gamma=0.5$ (blue), $\beta=0.2,\gamma=0.8$
      (green), and $\beta=0.6, \gamma=0.8$ (black). The right plot
      shows the same four curves
      with the number of SIS trials plotted on a logarithmic scale.  Note that
      the algorithm appears
      to be converging in sub-exponential time.
      Recall from Figure \ref{fig:convergence} that it is
     converging to the wrong estimate.}
       \label{fig:compareBad}

\end{figure}

Figure \ref{fig:compareReg} represents the number of trials
required by the SIS procedure (computed by our stopping heuristic)
on several examples for $n\times n$  matrices. The four curves
correspond to $5$, $10$, $\lfloor 5\log n\rfloor$ and $\lfloor
n/2\rfloor$-regular row and column sums. The $x$-axis represents
$n$, the number of rows and columns, and the $y$-axis captures the
required number of SIS trials. For each $n$ and each of these row
and column sums, we took $20$ independent runs and we plotted the
median number of trials. For comparison, in Figure
\ref{fig:compareBad} we plotted the estimated running time for our bad
example from Theorem \ref{thm:verybad} (recall that this
is likely the running time needed to converge to a wrong value!)
for $n+m$ ranging from $20$ to $140$ and various settings of $\beta,\gamma$:
$0.1,0.5$ (red), $0.5,0.5$ (blue), $0.2,0.8$ (green), and
$0.6,0.8$ (black).  In this case it is clear that the convergence time is
considerably slower
compared with the examples in Figure \ref{fig:compareReg}.

\subsection{Alternating Rows and Columns}
\label{sec:altRowCol}

The bad input instances from Theorem \ref{thm:verybad} can be
efficiently handled by an enhanced SIS approach which considers
both rows and columns for updating.  More precisely, the enhanced
SIS algorithm
assigns entries to the row {\it or\/} column with the largest residual sum.
We believe there are input instances for which this enhanced SIS
algorithm requires exponential time, but proving such a result appears
to be technically difficult.  We instead give experimental evidence that
there are such bad input instances for the enhanced SIS algorithm.

Specifically, we conjecture that for the family of inputs of the
form $\mathbf{r} = (1,1,\dots,1,\lfloor{m/2}\rfloor,$ 
$\lfloor
m/2\rfloor,\dots,\lfloor m/2\rfloor)$ and $\mathbf{c} =
(1,1,\dots,1,\lfloor{m/2}\rfloor,\lfloor m/2\rfloor,\dots,\lfloor
m/2\rfloor)$, where $m$ denotes the overall number of rows and there
are $\lfloor{m/2}\rfloor$ rows with sum $\lfloor{m/2}\rfloor$ and
$\lfloor{m/2}\rfloor$ columns with sum $\lfloor{m/2}\rfloor$, the
enhanced SIS strategy fails to converge quickly to $|\Omega|$.  A
theoretical analysis of the performance on this family of inputs is
difficult because, unlike the simpler instances of
Theorem~\ref{thm:verybad},  the true row and column distributions
are apparently rather hard to estimate in this case. Therefore, we
opted to perform experiments that suggest that even the enhanced SIS
algorithm is inefficient for this class of inputs.

\begin{figure}[H]
    \begin{center}
      \includegraphics[type=eps,ext=.eps,read=.eps,scale=1]{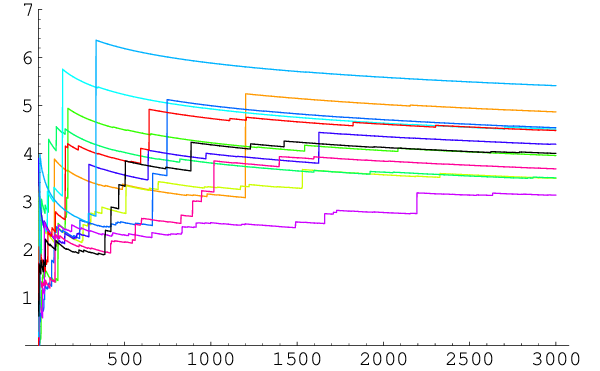}
    \end{center}
\caption{
  The estimate produced by the enhanced SIS scheme discussed in
Section \ref{sec:altRowCol} as a function of the number of trials
for the input $\mathbf{r}=\mathbf{c}=(1,1,\dots,1,\lfloor
m/2\rfloor,\lfloor m/2\rfloor,\dots,\lfloor m/2\rfloor)$, where $m$
is the number of rows (and columns), and the number of rows (and
columns) with marginal sum $\lfloor m/2 \rfloor$ is $\lfloor
m/2\rfloor$. The $x$-axis depicts the number of SIS trials scaled
down by a factor of 10,000, totalling 30,000,000 trials. The
$y$-axis depicts the SIS estimate divided by a factor of $10^{200}$,
on a logarithmic scale (base 10). Twelve independent runs are
shown.} \label{fig:altSIS}
\end{figure}

We now describe these experiments.  We did 30 million SIS trials for $m=100$,
and repeated this 12 times.  The estimates of $|\Omega|$ from these 12
experiments are presented in Figure~\ref{fig:altSIS}.  In the figure it
is clear that after 30 million trials these 12 experiments yield
quite different estimates of $|\Omega|$, differing by a factor
on the order of $10^3$.  These results strongly suggest that the enhanced SIS
algorithm has failed to converge to an estimate of $|\Omega|$ after 30 million
trials.  Moreover, we believe that the values produced by the
enhanced SIS algorithm after 30 million trials are substantial underestimates of
the true value of $|\Omega|$; however, since we know of no feasible method
for accurately estimating~$|\Omega|$ on these examples, we cannot compare
the experimental estimates to the true value of $|\Omega|$.

\begin{thebibliography}{99}
\bibitem{BKS}
M. Bayati, J.-H. Kim, and A. Saberi. A Sequential Algorithm for
Generating Random Graphs. {\em Algorithmica}, 58(4):860--910,
2010.

\bibitem{BesagClifford}
J. Besag and P. Clifford. Sequential Monte Carlo $p$-values. {\em
Biometrika}, 78(2):301--304, 1991.

\bibitem{BBV}
I. Bez\'akov\'a, N. Bhatnagar, and E. Vigoda.
Sampling Binary Contingency Tables with a Greedy Start.
{\em Random Structures and Algorithms}, 30(1-2):168--205, 2007.

\bibitem{arxiv-version}
I. Bez\'akov\'a,
A. Sinclair,
D. \v{S}tefankovi\v{c},
and E. Vigoda.
Negative Examples for Sequential Importance Sampling of Binary Contingency Tables.
Version available on the arXiv at:
\href{http://arxiv.org/abs/math/0606650}{http://arxiv.org/abs/math/0606650}

\bibitem{Blanchet}
J. Blanchet.
Efficient Importance Sampling for Binary Contingency Tables.
{\em  Annals of Applied Probability}, 19(3):949--982, 2009.

\bibitem{Blitzstein}
Joseph Blitzstein and Persi Diaconis. A Sequential Importance
Sampling Algorithm for Generating Random Graphs with Prescribed
Degrees. {\em Internet Mathematics}, 6(4):489--522, 2010.

\bibitem{CDHL}
Y. Chen, P. Diaconis, S. Holmes, and J.S. Liu.
Sequential Monte Carlo Methods for Statistical Analysis of Tables.
{\em Journal of the American Statistical Association}, 100:109--120, 2005.

\bibitem{CDS}
Y. Chen, I. Dinwoodie, and S. Sullivant.
Sequential Importance Sampling For Multiway Tables.
{\em The Annals of Statistics}, 34(1):523--545, 2006.

\bibitem{DP09}
D. P. Dubhashi and A. Panconesi,
{\em Concentration of measure for the analysis of randomized algorithms.}
Cambridge University Press, New York, 2009.

\bibitem{IGLRousset2005}
M. De Iorio, R. C. Griffiths, R. Lebois, and F. Rousset.
Stepwise Mutation Likelihood Computation by Sequential Importance Sampling
in Subdivided Population Models. {\em Theoretical Population Biology},
68:41--53, 2005.

\bibitem{JSV}
M. Jerrum, A. Sinclair and E. Vigoda.
A Polynomial-time Approximation Algorithm for the Permanent of a
Matrix with Non-negative Entries. {\em Journal of the Association
for Computing Machinery}, 51(4):671--697, 2004.

\bibitem{SigProc}
J. Miguez, and P. M. Djuric.
Blind Equalization by Sequential
Importance Sampling. {\em Proceedings of the IEEE International
Symposium on Circuits and Systems}, 845--848, 2002.

\bibitem{Tutte}
W. T. Tutte
A short proof of the factor theorem for finite graphs.
{\em Canad. J. Math}, 6:347--352, 1954.

\bibitem{ZhangLiu2002}
J. L. Zhang, and J. S. Liu.
 A New Sequential Importance Sampling Method and its Application to the
 Two-dimensional Hydrophobic-Hydrophilic Model.
 {\em Journal of Chemical Physics}, 117(7):3492--3498, 2002.
\end{thebibliography}
\end{document}